# A GENERALIZATION OF THE LINDEBERG PRINCIPLE[1]

By Sourav Chatterjee

*University of California, Berkeley*

We generalize Lindeberg's proof of the central limit theorem to an invariance principle for arbitrary smooth functions of independent and weakly dependent random variables. The result is applied to get a similar theorem for smooth functions of exchangeable random variables. This theorem allows us to identify, for the first time, the limiting spectral distributions of Wigner matrices with exchangeable entries.

**1. Introduction and results.** J. W. Lindeberg's elegant proof of the central limit theorem [12], despite being in the shadow of Fourier analytic methods for a long time, is now well known. It was revived by Trotter [20] and has since been successfully used to derive CLT's in infinite-dimensional spaces, where the Fourier analytic methods are not so useful.

While the original Lindeberg method and its extensions compare the distributions of convolutions in great generality (the history of which is irrelevant to our discussion, so we refer to [15] for details), it soon becomes clear that the same principle works not only for sums, but for more general smooth functions as well. Comparison of $f(X_1, \ldots, X_n)$ and $f(Y_1, \ldots, Y_n)$ for polynomial $f$ has been examined by Rotar [16] and Mossel, O'Donnell and Oleszkiewicz [13], and for general smooth $f$ with bounded derivatives by Chatterjee [5]. In [5], it is shown how to apply the method to establish universality in physical models, including the Sherrington–Kirkpatrick model of spin glasses. It was recently observed by Toufic Suidan [18] that the results in [5] can be used to give an immediate proof of the universality of last passage percolation in thin rectangles (originally a result of [3] and [4]). Developed independently, the Mossel, O'Donnell and Oleszkiewicz paper [13] is another repository of very striking modern applications of this very old idea.

Received August 2005; revised March 2006.
[1]Supported in part by NSF Grant DMS-04-06042.

A closer examination of Lindeberg's method reveals that there is a direct generalization by which the independence of the coordinates can be dispensed with. The argument, which may be possible to guess once the theorem has been stated, will be given in Section 2.

THEOREM 1.1. *Suppose $\mathbf{X}$ and $\mathbf{Y}$ are random vectors in $\mathbb{R}^n$ with $\mathbf{Y}$ having independent components. For $1 \le i \le n$, let*

$$A_i := \mathbb{E}|\mathbb{E}(X_i|X_1,\ldots,X_{i-1}) - \mathbb{E}(Y_i)|,$$
$$B_i := \mathbb{E}|\mathbb{E}(X_i^2|X_1,\ldots,X_{i-1}) - \mathbb{E}(Y_i^2)|.$$

*Let $M_3$ be a bound on $\max_i(\mathbb{E}|X_i|^3 + \mathbb{E}|Y_i|^3)$. Suppose $f:\mathbb{R}^n \to \mathbb{R}$ is a thrice continuously differentiable function, and for $r = 1,2,3$, let $L_r(f)$ be a finite constant such that $|\partial_i^r f(\mathbf{x})| \le L_r(f)$ for each $i$ and $\mathbf{x}$, where $\partial_i^r$ denotes the $r$-fold derivative in the $i$th coordinate. Then*

$$|\mathbb{E}f(\mathbf{X}) - \mathbb{E}f(\mathbf{Y})| \le \sum_{i=1}^n (A_i L_1(f) + \tfrac{1}{2} B_i L_2(f)) + \tfrac{1}{6} n L_3(f) M_3.$$

Let us now say a bit about the condition of boundedness of third derivatives. The implications of this condition have been inspected in detail in the context of convolutions by Zolotarev (see, e.g., [23]) and other authors. Zolotarev defines the $\zeta_3$-metric on the space of distributions as follows: $\zeta_3(F,G) = \sup_f |\int f\, dF - \int f\, dG|$, where the sup is taken over all $f$ with third derivative bounded by 1. The $\zeta_3$ metric has not been so popular in practice because of the difficulty in connecting this metric with the common notions of distance between measures.

However, instead of taking supremum over a class of $f$'s, we consider only individual functions of interest. For instance, in the random matrix scenario, our $f$ will be the Stieltjes transform of a matrix at a fixed $z \in \mathbb{C}\backslash\mathbb{R}$, which is a nice $C^\infty$ function of the original matrix. In the paper [5], the author considered the partition function of a disordered physical system (the Sherrington–Kirkpatrick model of spin glasses), which again turns out to be a $C^\infty$ function of the disorder matrix, and has nicely bounded derivatives.

The condition of boundedness of the derivatives can be dropped (as demonstrated in [16] and [13]) by careful examination of the remainder term; but our focus is different: We are more concerned with ways to extend the method to the case of weakly dependent variables (as in Theorem 1.1 above), and more specifically, to exchangeable random variables, as below.

*Exchangeable random variables.* Let us now present a surprisingly nontrivial application of the basic tool developed in the previous section. Suppose $\mathbf{X}$ is a vector with exchangeable components. Certainly, we cannot expect to replace the components of $\mathbf{X}$ by independent Gaussians as we did in Theorem 1.1. For instance, all the components may be equal to the same random variable, in which case there is no hope of replacing these variables by something generic. However, not all is lost; our next theorem shows that the following "summarization" of $\mathbf{X}$ can still be carried through:

Suppose $\mathbf{X}$ is a vector with exchangeable components, having finite fourth moments. Let

$$(1) \qquad \hat{\mu} := \frac{1}{n}\sum_{i=1}^{n} X_i \quad \text{and} \quad \hat{\sigma}^2 := \frac{1}{n}\sum_{i=1}^{n}(X_i - \hat{\mu})^2.$$

Let $\mathbf{Z}$ be a standard Gaussian vector in $\mathbb{R}^n$, independent of $\mathbf{X}$. Let $\bar{Z} := \frac{1}{n}\sum_{i=1}^{n} Z_i$ and

$$(2) \qquad Y_i := \hat{\mu} + \hat{\sigma}(Z_i - \bar{Z}), \qquad i = 1, \ldots, n.$$

Then, for sufficiently well-behaved $f$ (to be described below), we have $\mathbb{E}f(\mathbf{X}) \approx \mathbb{E}f(\mathbf{Y})$. That is, $\mathbf{X}$ can be "replaced" by the modified vector $\mathbf{Y}$ for evaluation under suitably smooth $f$. Note that in the process, we summarized the random vector $\mathbf{X}$ into the couple $(\hat{\mu}, \hat{\sigma})$. The precise statement is as follows:

THEOREM 1.2. *Suppose $\mathbf{X}$ is a random vector with exchangeable components, and $\hat{\mu}$, $\hat{\sigma}$ and $\mathbf{Y}$ are defined as in (1) and (2). Let $f : \mathbb{R}^n \to \mathbb{R}$ be a thrice continuously differentiable function, and for $r = 1, 2, 3$, let $L'_r(f)$ be a uniform bound on all $r$th partial derivatives of $f$, including mixed partials. For each $p$, let $m_p = \mathbb{E}|X_1 - \hat{\mu}|^p$. Then we have the bound*

$$|\mathbb{E}f(\mathbf{X}) - \mathbb{E}f(\mathbf{Y})| \leq 9.5 m_4^{1/2} L'_2(f) n^{1/2} + 13 m_3 L'_3(f) n.$$

We postpone the (somewhat long) proof of this theorem until Section 3, giving only a brief sketch at this point. The first step is to show that there is no loss of generality in assuming that $\hat{\mu} \equiv 0$ and $\hat{\sigma} \equiv 1$. Having assumed that, if we define

$$R_i = X_i + \frac{1}{n-i+1}\sum_{j=1}^{i-1} X_j,$$

then it is a straightforward exercise (which we will work out, nevertheless) that $\mathbb{E}(R_i | \mathcal{F}_{i-1}) = 0$, where $\mathcal{F}_{i-1}$ is the sigma-algebra generated by $X_1, \ldots, X_{i-1}$. The next step is to prove

$$\mathbb{E}(R_i^2 | \mathcal{F}_{i-1}) = 1 + O((n-i+1)^{-1/2}),$$

which is computationally slightly harder. Having established these approximations, we can now replace $R_i$'s by independent Gaussian variables $V_1, \ldots, V_n$, using Theorem 1.1. However, the inverse transform, which takes $\mathbf{R}$ to $\mathbf{X}$, does not take $\mathbf{V}$ to $\mathbf{Y}$ or anything resembling $\mathbf{Y}$, but to something that is close to $\mathbf{Y}$ in distribution. This will be formalized using Gaussian interpolation techniques.

Before moving on to applications, let us quickly mention without detailed justification that Theorem 1.2 has no apparent connection with de Finetti's theorem [6] or its finite version due to Diaconis and Freedman [8].

*An application*: *Wigner's law for exchangeable random variables*. Let us begin with a very quick introduction to some necessary material from random matrix theory.

*Spectral measures.* The *empirical spectral distribution* (ESD) of an $N \times N$ square matrix $A$ is the probability distribution $\frac{1}{N}\sum_{i=1}^{N}\delta_{\lambda_i}$, where $\lambda_1, \ldots, \lambda_N$ are the eigenvalues of $A$ repeated by multiplicities, and $\delta_x$ denotes the point mass at $x$. The weak limit of a sequence of ESDs is called the *limiting spectral distribution* (LSD) of the corresponding sequence of matrices. The existence and identification of LSDs for various kinds of random matrices is one of the main goals of random matrix theory. For a Hermitian matrix, the ESD is supported on the real line and hence has a corresponding cumulative distribution function. We will denote the c.d.f. for the ESD of a Hermitian matrix $A$ by $F_A$. Explicitly, $F_A(x) = \frac{1}{N}\#\{i : \lambda_i \le x\}$.

*Wigner matrices.* A standard Gaussian Wigner matrix of order $N$ is a matrix of the form $A_N = (N^{-1/2}X_{ij})_{1 \le i,j \le N}$, where $(X_{ij})_{1 \le i \le j \le N}$ is a collection of i.i.d. standard Gaussian random variables, and $X_{ij} = X_{ji}$ for $i > j$. Wigner [21] showed that the LSD for a sequence of standard Gaussian Wigner matrices (with order $N \to \infty$) is the semicircle law, which has density $(2\pi)^{-1}\sqrt{4-x^2}$ in the interval $[-2, 2]$.

It was later shown that the distribution of the entries does not play a significant role; convergence to the semicircle law holds under more general conditions (cf. [1, 2, 10]). The weakest known condition was given by Pastur [14]. It is claimed that the condition was shown to be necessary by Girko [9]. Although most conditions require independence of the entries on and above the diagonal, there have been some advances (e.g., [7, 17]) allowing certain kinds of dependence. However, none of these cover the case of exchangeable entries.

For a detailed exposition of the key results about the spectra of Wigner matrices and other results in the study of the spectral behavior of large random matrices, see [2] or [11].

Here we consider the question of identifying the limiting spectral distributions for Wigner matrices with exchangeable entries. The following theorem gives a precise answer under minimal assumptions.

THEOREM 1.3. *Suppose that for each $N$ we have a random matrix $A_N = (N^{-1/2} X_{ij}^N)_{1 \leq i,j \leq N}$, where the collection $(X_{ij}^N)_{1 \leq i \leq j \leq N}$ is exchangeable and $X_{ij}^N = X_{ji}^N$ for $i > j$. Let*

$$\hat{\mu}_N := \frac{2}{N(N+1)} \sum_{i \leq j \leq N} X_{ij}^N \quad \text{and} \quad \hat{\sigma}_N^2 := \frac{2}{N(N+1)} \sum_{i \leq j \leq N} (X_{ij}^N - \hat{\mu}_N)^2.$$

*Assume that $\hat{\sigma}_N > 0$ a.s. for all $N$ and $\sup_{N \geq 1} \mathbb{E}|(X_{11}^N - \hat{\mu}_N)/\hat{\sigma}_N|^4 < \infty$. Then the empirical spectral distribution of $\hat{\sigma}_N^{-1} A_N$ converges weakly to the semicircle law in probability.*

We will derive the above result from a quantitative bound (Lemma 4.1) on the difference between Stieltjes transforms (to be discussed in Section 4). The proof will show that it is actually enough to assume the weaker condition that $\mathbb{E}|(X_{11}^N - \hat{\mu}_N)/\hat{\sigma}_N|^4 = o(N^{2/3})$ as $N \to \infty$. It will also be evident that the argument can be adapted to more complicated exchangeability assumptions than the most basic one assumed above.

**2. Proof of Theorem 1.1.** Throughout the remainder of this article, we will use the notation $\partial_i f$ instead of the more familiar $\frac{\partial f}{\partial x_i}$. Similarly, we will write $\partial_i \partial_j f$ instead of $\frac{\partial^2 f}{\partial x_i \partial x_j}$ and so on.

Now let us begin with the proof. Without loss of generality, we can assume that $\mathbf{X}$ and $\mathbf{Y}$ are defined on the same probability space and are independent. For each $i$, $0 \leq i \leq n$, let

$$\mathbf{Z}_i = (X_1, \ldots, X_i, Y_{i+1}, \ldots, Y_n) \quad \text{and} \quad \mathbf{Z}_i^0 = (X_1, \ldots, X_{i-1}, 0, Y_{i+1}, \ldots, Y_n).$$

Then clearly

$$\mathbb{E} f(\mathbf{X}) - \mathbb{E} f(\mathbf{Y}) = \sum_{i=1}^n (\mathbb{E} f(\mathbf{Z}_i) - \mathbb{E} f(\mathbf{Z}_{i-1})).$$

Now, by third-order Taylor approximation,

$$\left| f(\mathbf{Z}_i) - f(\mathbf{Z}_i^0) - X_i \, \partial_i f(\mathbf{Z}_i^0) - \frac{X_i^2}{2} \partial_i^2 f(\mathbf{Z}_i^0) \right| \leq \frac{|X_i|^3 L_3(f)}{6}$$

and similarly,

$$\left| f(\mathbf{Z}_{i-1}) - f(\mathbf{Z}_i^0) - Y_i \, \partial_i f(\mathbf{Z}_i^0) - \frac{Y_i^2}{2} \partial_i^2 f(\mathbf{Z}_i^0) \right| \leq \frac{|Y_i|^3 L_3(f)}{6}.$$

Now, since the $Y_i$'s are independent, we have

$$\mathbb{E}((X_i - Y_i) \, \partial_i f(\mathbf{Z}_i^0)) = \mathbb{E}((\mathbb{E}(X_i | X_1, \ldots, X_{i-1}) - \mathbb{E}(Y_i)) \, \partial_i f(\mathbf{Z}_i^0)).$$

Similarly, we also have
$$\mathbb{E}((X_i^2 - Y_i^2)\,\partial_i^2 f(\mathbf{Z}_i^0)) = \mathbb{E}((\mathbb{E}(X_i^2|X_1,\ldots,X_{i-1}) - \mathbb{E}(Y_i^2))\,\partial_i^2 f(\mathbf{Z}_i^0)).$$

Thus, for any $i$,
$$\begin{aligned}|\mathbb{E}f(\mathbf{Z}_i) - \mathbb{E}(\mathbf{Z}_{i-1})| &\leq \tfrac{1}{6}L_3(f)(\mathbb{E}|X_i|^3 + \mathbb{E}|Y_i|^3) \\ &\quad + |\mathbb{E}((X_i - Y_i)\,\partial_i f(\mathbf{Z}_i^0))| + \tfrac{1}{2}|\mathbb{E}((X_i^2 - Y_i^2)\partial_i^2 f(\mathbf{Z}_i^0))| \\ &\leq \tfrac{1}{6}L_3(f)M_3 + A_i L_1(f) + \tfrac{1}{2}B_i L_2(f).\end{aligned}$$

This completes the proof.

**3. Proof of Theorem 1.2.** First, note that $\mathbf{X} = \mathbf{Y}$ on the event $\{\hat{\sigma} = 0\}$. Thus, if $\mathbb{P}\{\hat{\sigma} = 0\} = 1$, there is nothing to prove. If $\mathbb{P}\{\hat{\sigma} = 0\} < 1$, we can condition on the event $\{\hat{\sigma} > 0\}$ and consequently assume, without loss of generality, that $\hat{\sigma} > 0$ almost surely, because the conditioning retains the exchangeability of the $X_i$'s. Thus, let us assume that $\mathbb{P}\{\hat{\sigma} > 0\} = 1$.

For $i = 1, \ldots, n$ let $\tilde{X}_i = (X_i - \hat{\mu})/\hat{\sigma}$. Then $\sum_{i=1}^n \tilde{X}_i = 0$ and $\sum_{i=1}^n \tilde{X}_i^2 = n$ (we will be using these identities numerous times, often without mention). Let $\tilde{Z}_i = Z_i - \bar{Z}$, where $\bar{Z} = \frac{1}{n}\sum_{i=1}^n Z_i$.

In the following, we will use $\mathbb{E}_0$ and $\mathbb{P}_0$ to denote the expectation and probability conditional on the pair $(\hat{\mu}, \hat{\sigma})$. Observe that $\tilde{\mathbf{X}}$ is a vector with exchangeable components under $\mathbb{P}_0$ for all values of $(\hat{\mu}, \hat{\sigma})$.

Now assume that $(\hat{\mu}, \hat{\sigma})$ is given and fixed. Let
$$f_0(x_1, \ldots, x_n) := f(\hat{\mu} + \hat{\sigma}x_1, \ldots, \hat{\mu} + \hat{\sigma}x_n).$$
Then $f(\mathbf{X}) = f_0(\tilde{\mathbf{X}})$ and $f(\mathbf{Y}) = f_0(\tilde{\mathbf{Z}})$. Note that $L'_r(f_0) = L'_r(f)\hat{\sigma}^r$, where $L'_r(g)$ denotes a uniform bound on the $r$th order derivatives of a function $g$, including mixed partials.

First, we need to do a list of computations. For $0 \leq i \leq n$, let $\mathcal{F}_i$ be the sigma-algebra generated by $\{\tilde{X}_1, \ldots, \tilde{X}_i\}$ and $(\hat{\mu}, \hat{\sigma})$. Since the $\tilde{X}_j$'s are exchangeable given $(\hat{\mu}, \hat{\sigma})$, we have $\mathbb{E}_0(\tilde{X}_k|\mathcal{F}_{i-1}) = \mathbb{E}_0(\tilde{X}_l|\mathcal{F}_{i-1})$ for every $k, l > i - 1$, and hence
$$\mathbb{E}_0(\tilde{X}_i|\mathcal{F}_{i-1}) = \mathbb{E}_0\left(\frac{1}{n-i+1}\sum_{j=i}^n \tilde{X}_j \bigg| \mathcal{F}_{i-1}\right).$$

Now, since $\sum_{j=1}^n \tilde{X}_j = 0$,
$$\frac{1}{n-i+1}\sum_{j=i}^n \tilde{X}_j = -\frac{1}{n-i+1}\sum_{j=1}^{i-1} \tilde{X}_j,$$
which is $\mathcal{F}_{i-1}$-measurable. Thus,

(3) $$\mathbb{E}_0(\tilde{X}_i|\mathcal{F}_{i-1}) = -\frac{1}{n-i+1}\sum_{j=1}^{i-1} \tilde{X}_j = \frac{1}{n-i+1}\sum_{j=i}^n \tilde{X}_j.$$

From the above identity and exchangeability, it follows that

$$\mathbb{E}_0(\mathbb{E}_0(\tilde{X}_i|\mathcal{F}_{i-1})^2) = \frac{1}{(n-i+1)^2}[(n-i+1)\mathbb{E}_0(\tilde{X}_1^2)$$
$$+ (n-i+1)(n-i)\mathbb{E}_0(\tilde{X}_1\tilde{X}_2)].$$

Now clearly $\mathbb{E}_0(\tilde{X}_1^2) = \mathbb{E}_0(\frac{1}{n}\sum_{i=1}^n \tilde{X}_i^2) = 1$ and

$$\mathbb{E}_0(\tilde{X}_1\tilde{X}_2) = \frac{1}{n(n-1)} \sum_{i \neq j} \mathbb{E}_0(\tilde{X}_i\tilde{X}_j)$$
$$= -\frac{1}{n(n-1)} \sum_{i=1}^n \mathbb{E}_0(\tilde{X}_i^2) = -\frac{1}{n-1}.$$

(The second equality holds because $\sum_{j=1}^n \tilde{X}_j = 0$.) Combining, we get

(4) $$\mathbb{E}_0(\mathbb{E}_0(\tilde{X}_i|\mathcal{F}_{i-1})^2) = \frac{(i-1)}{(n-i+1)(n-1)}.$$

Again, by a similar argument as before (using the identity $\sum_{i=1}^n \tilde{X}_i^2 = n$), we have

$$\mathbb{E}_0(\tilde{X}_i^2|\mathcal{F}_{i-1}) = \frac{1}{n-i+1} \sum_{j=i}^n \tilde{X}_j^2.$$

It follows that

$$\mathrm{Var}_0(\mathbb{E}_0(\tilde{X}_i^2|\mathcal{F}_{i-1})) = \frac{1}{(n-i+1)^2}[(n-i+1)\mathrm{Var}_0(\tilde{X}_1^2)$$
$$+ (n-i+1)(n-i)\mathrm{Cov}_0(\tilde{X}_1^2, \tilde{X}_2^2)].$$

Now, $\mathrm{Var}_0(\tilde{X}_1^2) \leq \mathbb{E}_0(\tilde{X}_1^4)$, and

$$\mathrm{Cov}_0(\tilde{X}_1^2, \tilde{X}_2^2) = \frac{1}{n(n-1)} \sum_{i \neq j} \mathbb{E}_0(\tilde{X}_i^2 \tilde{X}_j^2) - \mathbb{E}_0(\tilde{X}_1^2)\mathbb{E}_0(\tilde{X}_2^2)$$
$$= \frac{1}{n(n-1)} \sum_{i=1}^n \mathbb{E}_0(\tilde{X}_i^2(n - \tilde{X}_i^2)) - 1$$
$$= \frac{1 - \mathbb{E}_0(\tilde{X}_1^4)}{n-1} \leq 0, \qquad \text{since } \mathbb{E}_0(\tilde{X}_1^4) \geq (\mathbb{E}_0(\tilde{X}_1^2))^2 = 1.$$

Combining, we get

(5) $$\mathrm{Var}_0(\mathbb{E}_0(\tilde{X}_i^2|\mathcal{F}_{i-1})) \leq \frac{\mathbb{E}_0(\tilde{X}_1^4)}{n-i+1}.$$

Now let $G$ be the matrix whose $(i,j)$th element is

$$g_{ij} = \begin{cases} 1/(n-i+1), & \text{if } i > j, \\ 1, & \text{if } i = j, \\ 0, & \text{if } i < j. \end{cases}$$

Let $\mathbf{R} = G\tilde{\mathbf{X}}$. Then $R_i$ is $\mathcal{F}_i$-measurable, and from (3), we see that

(6) $$\mathbb{E}_0(R_i|\mathcal{F}_{i-1}) = 0.$$

Next, note that

$$\mathbb{E}_0(R_i^2|\mathcal{F}_{i-1}) = \mathbb{E}_0((\tilde{X}_i - \mathbb{E}_0(\tilde{X}_i|\mathcal{F}_{i-1}))^2|\mathcal{F}_{i-1})$$
$$= \mathbb{E}_0(\tilde{X}_i^2|\mathcal{F}_{i-1}) - (\mathbb{E}_0(\tilde{X}_i|\mathcal{F}_{i-1}))^2.$$

Using (4), (5), the triangle inequality and the fact that $\mathbb{E}_0(\tilde{X}_i^2) = 1$, we get

$$\mathbb{E}_0|\mathbb{E}_0(R_i^2|\mathcal{F}_{i-1}) - 1| \leq \mathbb{E}_0|\mathbb{E}_0(\tilde{X}_i^2|\mathcal{F}_{i-1}) - 1| + \mathbb{E}_0(\mathbb{E}_0(\tilde{X}_i|\mathcal{F}_{i-1})^2)$$

(7) $$\leq \sqrt{\frac{\mathbb{E}_0(\tilde{X}_1^4)}{n-i+1}} + \frac{(i-1)}{(n-i+1)(n-1)}$$

$$\leq 2\sqrt{\frac{\mathbb{E}_0(\tilde{X}_1^4)}{n-i+1}} \quad \text{since } \mathbb{E}_0(\tilde{X}_1^4) \geq 1.$$

We will now temporarily use the notation $\|R_i\|_3^0$ for $(\mathbb{E}_0|R_i|^3)^{1/3}$, which is the conditional $L^3$ norm of $R_i$. By Minkowski's inequality and exchangeability, we have

$$\|R_i\|_3^0 \leq \|\tilde{X}_i\|_3^0 + \frac{1}{n-i+1}\sum_{j=i}^n \|\tilde{X}_j\|_3^0 = 2\|\tilde{X}_1\|_3^0.$$

This bound can be rewritten as

(8) $$\mathbb{E}_0|R_i|^3 \leq 8\hat{\sigma}^{-3}\mathbb{E}_0|X_1 - \hat{\mu}|^3.$$

Now define the function $f_1 : \mathbb{R}^n \to \mathbb{R}$ as $f_1(\mathbf{x}) := f_0(G^{-1}\mathbf{x})$. Then $f_1(\mathbf{R}) = f_0(\tilde{\mathbf{X}})$. Let $g^{ij}$ denote the $(i,j)$th element of $G^{-1}$. It is a simple exercise to verify that

$$g^{ij} = \begin{cases} -1/(n-j), & \text{if } i > j, \\ 1, & \text{if } i = j, \\ 0, & \text{if } i < j. \end{cases}$$

Using the chain rule we see that for any $j$, $r$ and $\mathbf{x}$,

$$\partial_j^r f_1(\mathbf{x}) = \sum_{1 \leq i_1, \ldots, i_r \leq n} (\partial_{i_1}\partial_{i_2}\cdots\partial_{i_r} f_0)(G^{-1}\mathbf{x})g^{i_1 j}g^{i_2 j}\cdots g^{i_r j}.$$

Thus for any $r$,

(9) $$L_r(f_1) \leq L'_r(f_0) \max_{1 \leq j \leq n} \left( \sum_{i=1}^n |g^{ij}| \right)^r = L'_r(f_0) 2^r = L'_r(f) \hat{\sigma}^r 2^r,$$

where $L_r$ is defined as in the statement of Theorem 1.1. Now let $\mathbf{V}$ be a standard Gaussian vector in $\mathbb{R}^n$. Using the bounds from (6)–(9) in Theorem 1.1, we get

$$|\mathbb{E}_0 f_1(\mathbf{R}) - \mathbb{E}_0 f_1(\mathbf{V})|$$
$$\leq \frac{1}{2} L_2(f_1) \sum_{i=1}^n \mathbb{E}_0 |\mathbb{E}_0(R_i^2 | \mathcal{F}_{i-1}) - 1| + \frac{1}{6} n L_3(f_1) \max_{1 \leq i \leq n} (\mathbb{E}_0 |R_i|^3 + \mathbb{E}|V_i|^3)$$
$$\leq 4 L'_2(f) \hat{\sigma}^2 \sum_{i=1}^n \sqrt{\frac{\mathbb{E}_0(\tilde{X}_1^4)}{n-i+1}} + \frac{8}{6} n L'_3(f) \hat{\sigma}^3 (8 \hat{\sigma}^{-3} \mathbb{E}_0 |X_1 - \hat{\mu}|^3 + \mathbb{E}|V_1|^3).$$

Now $\mathbb{E}|\tilde{X}_1|^4 = \hat{\sigma}^{-4} \mathbb{E}_0 |X_1 - \hat{\mu}|^4$, and by comparing sums with integrals, we have $\sum_{i=1}^n (n-i+1)^{-1/2} \leq 2\sqrt{n}$. Thus,

(10) $$|\mathbb{E}_0 f_1(\mathbf{R}) - \mathbb{E}_0 f_1(\mathbf{V})|$$
$$\leq 8 L'_2(f) \sqrt{n \mathbb{E}_0 |X_1 - \hat{\mu}|^4}$$
$$+ \tfrac{4}{3} n L'_3(f) (8 \mathbb{E}_0 |X_1 - \hat{\mu}|^3 + \hat{\sigma}^3 \mathbb{E}|V_1|^3).$$

Let $\mathbf{U} = G^{-1} \mathbf{V}$. Explicitly,

$$U_i = V_i - \sum_{j=1}^{i-1} \frac{V_j}{n-j}.$$

Now note that $f_1(\mathbf{V}) = f_0(\mathbf{U})$, $f_0(\tilde{\mathbf{Z}}) = f(\mathbf{Y})$ and $f_1(\mathbf{R}) = f_0(\tilde{\mathbf{X}}) = f(\mathbf{X})$. Combining, we get

(11) $$|\mathbb{E}_0 f(\mathbf{X}) - \mathbb{E}_0 f(\mathbf{Y})| = |\mathbb{E}_0 f_1(\mathbf{R}) - \mathbb{E}_0 f_0(\tilde{\mathbf{Z}})|$$
$$\leq |\mathbb{E}_0 f_1(\mathbf{R}) - \mathbb{E}_0 f_1(\mathbf{V})|$$
$$+ |\mathbb{E}_0 f_0(\mathbf{U}) - \mathbb{E}_0 f_0(\tilde{\mathbf{Z}})|.$$

We already have a bound on $|\mathbb{E}_0 f_1(\mathbf{R}) - \mathbb{E}_0 f_1(\mathbf{V})|$ from (10). We will now compute a bound on $|\mathbb{E}_0 f_0(\mathbf{U}) - \mathbb{E}_0 f_0(\tilde{\mathbf{Z}})|$, where recall that $\tilde{Z}_i = Z_i - \bar{Z}$, and $\mathbf{Z}$ is a standard Gaussian vector. To do that, we first need to do some computations. Let $\tilde{\sigma}_{ij} := \text{Cov}(U_i, U_j)$. Then

$$\tilde{\sigma}_{ij} = \begin{cases} -(n-j)^{-1} + \sum_{k=1}^{j-1} (n-k)^{-2}, & \text{if } i > j, \\ 1 + \sum_{k=1}^{j-1} (n-k)^{-2}, & \text{if } i = j, \\ \tilde{\sigma}_{ji}, & \text{if } i < j. \end{cases}$$

Now, for $i > j$, we can rewrite the first term in $\tilde{\sigma}_{ij}$ as a telescoping sum to get

$$\tilde{\sigma}_{ij} = -\frac{1}{n-1} - \sum_{k=1}^{j-1}\left(\frac{1}{n-k-1} - \frac{1}{n-k}\right) + \sum_{k=1}^{j-1}\frac{1}{(n-k)^2}$$

$$= -\frac{1}{n-1} - \sum_{k=1}^{j-1}\frac{1}{(n-k)^2(n-k-1)}.$$

Thus, if we define

$$\sigma_{ij} := \operatorname{Cov}(\tilde{Z}_i, \tilde{Z}_j) = \begin{cases} -1/n, & \text{if } i \neq j, \\ (n-1)/n, & \text{if } i = j, \end{cases}$$

then

$$\sum_{i,j}|\sigma_{ij} - \tilde{\sigma}_{ij}| = \sum_{i=1}^{n}|\sigma_{ii} - \tilde{\sigma}_{ii}| + 2\sum_{i=1}^{n}\sum_{j=1}^{i-1}|\sigma_{ij} - \tilde{\sigma}_{ij}|$$

$$\leq 2 + \sum_{i=1}^{n}\sum_{k=1}^{i-1}\frac{1}{(n-k)^2}$$

(12)

$$+ 2\sum_{i=1}^{n}\sum_{j=1}^{i-1}\sum_{k=1}^{j-1}\frac{1}{(n-k)^2(n-k-1)}$$

$$= 2 + \sum_{k=1}^{n-1}\frac{1}{n-k} + \sum_{k=1}^{n-2}\frac{1}{n-k} = 3 + 2\sum_{k=2}^{n-1}\frac{1}{k}.$$

We will use the well-known "Gaussian interpolation technique" for bounding $|\mathbb{E}_0 f_0(\mathbf{U}) - \mathbb{E}_0 f_0(\tilde{\mathbf{Z}})|$. This classical method for proving Slepian-type inequalities has been used extensively in recent years by Talagrand [19] in his efforts to obtain a rigorous version of the cavity method for spin glasses. For each $t \in [0,1]$, let $\mathbf{W}_t = \sqrt{1-t}\mathbf{U} + \sqrt{t}\tilde{\mathbf{Z}}$. Then

$$\mathbb{E}_0 f_0(\tilde{\mathbf{Z}}) - \mathbb{E}_0 f_0(\mathbf{U})$$

(13)
$$= \mathbb{E}_0\left[\int_0^1 \frac{d}{dt} f_0(\mathbf{W}_t)\, dt\right]$$

$$= \mathbb{E}_0\left[\int_0^1 \sum_{i=1}^{n}\left(\frac{\tilde{Z}_i}{2\sqrt{t}} - \frac{U_i}{2\sqrt{1-t}}\right)\partial_i f_0(\mathbf{W}_t)\, dt\right].$$

Now, if a random vector $\boldsymbol{\xi} = (\xi_1, \ldots, \xi_n)$ has a centered Gaussian distribution, then it is not difficult to show using integration by parts that for any differentiable function $h$ with subexponential growth at infinity, and any $i$, the following identity holds:

$$\mathbb{E}(\xi_i h(\boldsymbol{\xi})) = \sum_{j=1}^n \mathbb{E}(\xi_i \xi_j)\mathbb{E}(\partial_j h(\boldsymbol{\xi})).$$

Since we do not want to expand our list of references, let us refer to Appendix A.6 of Talagrand's book [19] for a proof. Applying this result to our problem (after noting that interchanging integrals is not an issue since everything is bounded), we get

$$\mathbb{E}_0(U_i\,\partial_i f_0(\mathbf{W}_t)) = \sqrt{1-t}\sum_{j=1}^n \tilde{\sigma}_{ij}\mathbb{E}_0(\partial_j\,\partial_i f_0(\mathbf{W}_t))$$

and similarly,

$$\mathbb{E}_0(\tilde{Z}_i\,\partial_i f_0(\mathbf{W}_t)) = \sqrt{t}\sum_{j=1}^n \sigma_{ij}\mathbb{E}_0(\partial_j\,\partial_i f_0(\mathbf{W}_t)).$$

Combining, we have

$$\mathbb{E}_0 f_0(\tilde{\mathbf{Z}}) - \mathbb{E}_0 f_0(\mathbf{U}) = \tfrac{1}{2}\int_0^1 \sum_{1\le i,j\le n} \mathbb{E}_0[\partial_i\,\partial_j f_0(\mathbf{W}_t)](\sigma_{ij} - \tilde{\sigma}_{ij})\,dt.$$

Using the bound from (12), we get

(14) $$|\mathbb{E}_0 f_0(\mathbf{U}) - \mathbb{E}_0 f_0(\tilde{\mathbf{Z}})| \le \frac{1}{2}L_2'(f_0)\left(3 + 2\sum_{k=2}^{n-1}\frac{1}{k}\right).$$

Combining this with (10) and (14), we get

$$\begin{aligned}
|\mathbb{E}f(\mathbf{X}) - \mathbb{E}f(\mathbf{Y})| \\
\le \mathbb{E}|\mathbb{E}_0 f(\mathbf{X}) - \mathbb{E}_0 f(\mathbf{Y})| \\
\le 8L_2'(f)\sqrt{n\mathbb{E}|X_1 - \hat{\mu}|^4} \\
+ \frac{4}{3}nL_3'(f)(8\mathbb{E}|X_1 - \hat{\mu}|^3 + \mathbb{E}(\hat{\sigma}^3)\mathbb{E}|V_1|^3) \\
+ \frac{1}{2}L_2'(f)\mathbb{E}(\hat{\sigma}^2)\left(3 + 2\sum_{k=2}^{n-1}\frac{1}{k}\right).
\end{aligned}$$

To complete the proof, we apply Jensen's inequality to get $\mathbb{E}(\hat{\sigma}^r) \le \mathbb{E}|X_1 - \hat{\mu}|^r$ for $r \ge 2$ and use the crude bounds $3 + 2\sum_{k=2}^{n-1} k^{-1} \le 3\sqrt{n}$ and $\mathbb{E}|V_1|^3 \le 1.7$ to unify terms.

**4. Proof of Theorem 1.3.** We will now prove Theorem 1.3 via an application of Theorem 1.2, by using the Stieltjes transform of the spectral measure as a smooth function of the matrix entries. The Stieltjes transform (or Cauchy transform, or resolvent) of a cumulative probability distribution function $F$ on $\mathbb{R}$ is defined as

$$m_F(z) := \int_{-\infty}^{\infty} \frac{1}{x-z} dF(x) \quad \text{for every } z \in \mathbb{C}\backslash\mathbb{R}. \tag{15}$$

Analogously, the Stieltjes transform of an $N \times N$ Hermitian matrix $A$ at a number $z \in \mathbb{C}\backslash\mathbb{R}$ is defined as

$$m_A(z) := \frac{1}{N} \text{Tr}((A - zI)^{-1}), \tag{16}$$

where $I$ is the identity matrix of order $N$. Note that this is just the Stieltjes transform of the empirical spectral distribution (ESD) of $A$. The ESDs of a sequence $\{A_N\}_{N=1}^{\infty}$ of random Hermitian matrices converge weakly in probability to a distribution $F$ if and only if

$$m_{A_N}(z) \xrightarrow{P} m_F(z) \quad \text{for every } z \in \mathbb{C}\backslash\mathbb{R}.$$

For the proof of this result and further details like Berry–Esseen-type error bounds, we refer to [2], pages 639–640.

Now recall that if $A(x)$ is a matrix-valued differentiable function of a scalar $x$, and $G(x) := (A(x) - zI)^{-1}$, where $z \in \mathbb{C}\backslash\mathbb{R}$ and $I$ is the identity matrix, then

$$\frac{dG}{dx} = -G \frac{dA}{dx} G. \tag{17}$$

This standard result is obtained by differentiating both sides of the identity $G(A - zI) \equiv I$. Differentiability follows from the fact that the elements of the inverse of a matrix are all rational functions of the elements of the original matrix. Higher-order derivatives may be computed by repeatedly applying the above formula.

The following lemma is the key to the proof of Theorem 1.3:

LEMMA 4.1. *Suppose that for each $N$ we have a random matrix $\tilde{A}_N = (N^{-1/2} \tilde{X}_{ij}^N)_{1 \le i,j \le N}$, where the collection $(\tilde{X}_{ij}^N)_{1 \le i \le j \le N}$ is exchangeable and $\tilde{X}_{ij}^N = \tilde{X}_{ji}^N$ for $i > j$. Suppose*

$$\frac{2}{N(N+1)} \sum_{i \le j \le N} \tilde{X}_{ij}^N = 0 \quad \text{and} \quad \frac{2}{N(N+1)} \sum_{i \le j \le N} (\tilde{X}_{ij}^N)^2 = 1 \quad a.s.$$

*For each $N$, let $(Z_{ij}^N)_{1 \le i \le j \le N}$ be a collection of i.i.d. standard Gaussian random variables and let $Y_{ij}^N = Z_{ij}^N - \bar{Z}^N$, where $\bar{Z}^N = \frac{2}{N(N+1)} \sum_{i \le j} Z_{ij}^N$.*

Let $Y_{ij}^N = Y_{ji}^N$ for $i > j$, and let $B_N = (N^{-1/2} Y_{ij}^N)_{1 \leq i,j \leq N}$. Then, for any $z \in \mathbb{C} \backslash \mathbb{R}$, and any $g : \mathbb{R} \to \mathbb{R}$ with bounded derivatives up to the third order, we have

$$|\mathbb{E}g(\operatorname{Re} m_{\tilde{A}_N}(z)) - \mathbb{E}g(\operatorname{Re} m_{B_N}(z))|$$
$$\leq C_1 N^{-1} (\mathbb{E}|\tilde{X}_{12}^N|^4)^{1/2} + C_2 N^{-1/2} \mathbb{E}|\tilde{X}_{12}^N|^3,$$

where $m$ is the Stieltjes transform as defined in (16) and $C_1$ and $C_2$ are constants depending only on $g$ and $z$. The quantity $|\mathbb{E}g(\operatorname{Im} m_{\tilde{A}_N}(z)) - \mathbb{E}g(\operatorname{Im} m_{B_N}(z))|$ also admits the same upper bound.

To complete the proof of Theorem 1.3 using this lemma, we need the following fact about spectral distributions of Hermitian matrices:

LEMMA 4.2 (Quoted from [2], Lemma 2.2). *Let $A$ and $B$ be two $N \times N$ Hermitian matrices, with empirical distribution functions $F_A$ and $F_B$. Then*

$$\|F_A - F_B\|_\infty \leq \frac{1}{N} \operatorname{rank}(A - B).$$

This lemma is an easy consequence of the well-known interlacing inequalities for eigenvalues of Hermitian matrices. Let us now complete the proof of Theorem 1.3 by combining Lemma 4.1 and Lemma 4.2.

PROOF OF THEOREM 1.3. Let $\tilde{X}_{ij}^N = (X_{ij}^N - \hat{\mu}_N)/\hat{\sigma}_N$, and let $\tilde{A}_N = (N^{-1/2} \tilde{X}_{ij}^N)_{1 \leq i,j \leq N}$. Clearly, $\tilde{A}_N$ satisfies the hypotheses of Lemma 4.1. Thus, we have

$$|\mathbb{E}g(\operatorname{Re} m_{\tilde{A}_N}(z)) - \mathbb{E}g(\operatorname{Re} m_{B_N}(z))|$$
$$\leq C_1 N^{-1} (\mathbb{E}|\tilde{X}_{12}^N|^4)^{1/2} + C_2 N^{-1/2} \mathbb{E}|\tilde{X}_{12}^N|^3.$$

The same bound holds for $|\mathbb{E}g(\operatorname{Im} m_{\tilde{A}_N}(z)) - \mathbb{E}g(\operatorname{Im} m_{B_N}(z))|$ as well. The bound converges to zero if $\mathbb{E}|\tilde{X}_{12}^N|^4 = o(N^{2/3})$, and thus under that condition, $\tilde{A}_N$ and $B_N$ must have the same LSD. Finally, observe that by Lemma 4.2, the sequence $\{\hat{\sigma}_N^{-1} A_N\}$ has the same LSD as $\{\tilde{A}_N\}$, and $F_{B_N}$ converges weakly to the semicircle distribution in probability. This completes the proof. $\square$

PROOF OF LEMMA 4.1. To formalize things in a way that is suitable for our purpose, consider the map $A$ which "constructs" Wigner matrices of order $N$. Let $n = N(N+1)/2$ and write elements of $\mathbb{R}^n$ as $\mathbf{x} = (x_{ij})_{1 \leq i \leq j \leq N}$. For any $\mathbf{x} \in \mathbb{R}^n$, let $A(\mathbf{x})$ be the matrix defined as

(18) $$A(\mathbf{x})_{ij} := \begin{cases} N^{-1/2} x_{ij}, & \text{if } i \leq j, \\ N^{-1/2} x_{ji}, & \text{if } i > j. \end{cases}$$

Now let us fix $z = u + \sqrt{-1}v \in \mathbb{C}$, with $v \neq 0$. Let $G(\mathbf{x}) := (A(\mathbf{x}) - zI)^{-1}$, and define $h : \mathbb{R}^n \to \mathbb{R}$ as

$$h(\mathbf{x}) := N^{-1} \operatorname{Tr}(G(\mathbf{x})).$$

For any $\alpha \in \{(i,j)\}_{1 \leq i \leq j \leq n}$, we will write $\partial_\alpha h$ for $\partial h / \partial x_\alpha$ by our usual convention. From (17), it follows that for any $\alpha$,

$$\partial_\alpha h = -N^{-1} \operatorname{Tr}(G(\partial_\alpha A)G). \tag{19}$$

Now note that for any $\alpha, \beta \in \{(i,j)\}_{1 \leq i \leq j \leq n}$, we have $\partial_\beta \partial_\alpha A \equiv 0$. An easy computation involving repeated applications of (17) to the above expression for $\partial_\alpha h$ gives, for any $\alpha, \beta, \gamma \in \{(i,j)\}_{1 \leq i \leq j \leq N}$,

$$\partial_\beta \partial_\alpha h = N^{-1} \sum_{\{\beta', \alpha'\} = \{\beta, \alpha\}} \operatorname{Tr}(G(\partial_{\beta'} A)G(\partial_{\alpha'} A)G), \tag{20}$$

$$\partial_\gamma \partial_\beta \partial_\alpha h = -N^{-1} \sum_{\{\gamma', \beta', \alpha'\} = \{\gamma, \beta, \alpha\}} \operatorname{Tr}(G(\partial_{\gamma'} A)G(\partial_{\beta'} A)G(\partial_{\alpha'} A)G). \tag{21}$$

Note that the first sum runs over all permutations of $(\beta, \alpha)$, which amounts to only two terms. Similarly, the second sum involves six terms.

To bound (19), first note that $\operatorname{Tr}(G(\partial_\alpha A)G) = \operatorname{Tr}((\partial_\alpha A)G^2)$. Since $G^2$ has a spectral decomposition and all its eigenvalues are bounded by $|v|^{-2}$ in magnitude, it follows in particular that the elements of $G^2$ are also bounded by $|v|^{-2}$. Now, $\partial_\alpha A$ has at most two nonzero elements, which are equal to $N^{-1/2}$. Hence,

$$|\operatorname{Tr}(G(\partial_\alpha A)G)| = |\operatorname{Tr}((\partial_\alpha A)G^2)| \leq 2|v|^{-2} N^{-1/2}.$$

To bound (20) and (21), we need to recall the properties of the Hilbert–Schmidt norm for matrices. For an $N \times N$ complex matrix $B = (b_{ij})_{1 \leq i,j \leq N}$, the *Hilbert–Schmidt norm* of $B$ is defined as $\|B\| := (\sum_{i,j} |b_{ij}|^2)^{1/2}$. Besides the usual properties of a matrix norm, it has the following additional features: (a) $|\operatorname{Tr}(BC)| \leq \|B\|\|C\|$, (b) if $U$ is a unitary matrix, then for any $C$ of the same order, $\|CU\| = \|UC\| = \|C\|$, and (c) for a matrix $B$ admitting a spectral decomposition (i.e., a normal matrix) with eigenvalues $\lambda_1, \ldots, \lambda_N$, and any other matrix $C$ of the same order, $\max\{\|BC\|, \|CB\|\} \leq \max_{1 \leq i \leq N} |\lambda_i| \cdot \|C\|$. For a proof of these standard facts one can look up, for example, [22], pages 55–58.

Clearly, $G$ and the derivatives of $A$ are all normal matrices. Moreover, the eigenvalues of $G$ are bounded by $|v|^{-1}$, where $v = \operatorname{Im} z$. Thus, by the properties of the Hilbert–Schmidt norm listed above, we have

$$|\operatorname{Tr}(G(\partial_\beta A)G(\partial_\alpha A)G)| \leq \|G(\partial_\beta A)\| \|G(\partial_\alpha A)G\|$$
$$\leq |v|^{-3} \|\partial_\beta A\| \|\partial_\alpha A\|$$
$$\leq 2|v|^{-3} N^{-1}.$$

Similarly,

$$|\text{Tr}(G(\partial_\gamma A)G(\partial_\beta A)G(\partial_\alpha A)G)| \leq \|G(\partial_\gamma A)\|\|G(\partial_\beta A)G(\partial_\alpha A)G\|$$
$$\leq \|G(\partial_\gamma A)\|\|G(\partial_\beta A)G\|\|(\partial_\alpha A)G\|$$
$$\leq |v|^{-4}\|\partial_\gamma A\|\|\partial_\beta A\|\|\partial_\alpha A\|$$
$$\leq 2^{3/2}|v|^{-4}N^{-3/2}.$$

Finally, note that since the matrix entries are all real, therefore $\partial_\alpha(\text{Re}\,h) = \text{Re}\,\partial_\alpha h$ and so on. Thus, if we let $f = g \circ (\text{Re}\,h)$, then substituting the bounds obtained above in (19), (20) and (21), we get $L_2'(f) \leq K_1 N^{-2}$ and $L_3'(f) \leq K_2 N^{-5/2}$, where $K_1$ and $K_2$ are constants depending only on $g$ and $z$. By Theorem 1.2, it now follows that

$$|\mathbb{E}f(\tilde{\mathbf{X}}) - \mathbb{E}f(\mathbf{Y})| \leq 9.5 K_1 N^{-1}(\mathbb{E}|\tilde{X}_{12}|^4)^{1/2} + 13 K_2 N^{-1/2}\mathbb{E}|\tilde{X}_{12}|^3,$$

where $Y_{ij} = Z_{ij} - \bar{Z}$, and $Z_{ij}$'s are i.i.d. standard Gaussian random variables. This completes the proof of the lemma. □

**Acknowledgments.** The author thanks the anonymous referee for his very detailed and helpful comments on the material and the style of presentation. He is also grateful to his former thesis advisor Persi Diaconis for constant encouragement and support.

## REFERENCES

[1] ARNOLD, L. (1967). On the asymptotic distribution of the eigenvalues of random matrices. *J. Math. Anal. Appl.* **20** 262–268. MR0217833
[2] BAI, Z. D. (1999). Methodologies in the spectral analysis of large-dimensional random matrices, a review. *Statist. Sinica* **9** 611–677. MR1711663
[3] BAIK, J. and SUIDAN, T. M. (2005). A GUE central limit theorem and universality of directed first and last passage site percolation. *Int. Math. Res. Not.* **6** 325–337. MR2131383
[4] BODINEAU, T. and MARTIN, J. (2005). A universality property for last-passage percolation paths close to the axis. *Electron. Comm. Probab.* **10** 105–112. MR2150699
[5] CHATTERJEE, S. (2004). A simple invariance theorem. Available at http://arxiv.org/math.PR/0508213.
[6] DE FINETTI, B. (1969). Sulla prosequibilità di processi aleatori scambiabili. *Rend. Ist. Mat. Univ. Trieste* **1** 53–67. MR0292146
[7] BOUTET DE MONVEL, A. and KHORUNZHY, A. (1998). Limit theorems for random matrices. *Markov Process. Related Fields* **4** 175–197. MR1641617
[8] DIACONIS, P. and FREEDMAN, D. (1980). Finite exchangeable sequences. *Ann. Probab.* **8** 745–764. MR0577313
[9] GIRKO, V. L. (1988). *Spectral Theory of Random Matrices*. Nauka, Moscow. MR0955497
[10] GRENANDER, U. (1963). *Probabilities on Algebraic Structures*. Wiley, New York. MR0259969

[11] KHORUNZHY, A. M., KHORUZHENKO, B. A. and PASTUR, L. A. (1996). Asymptotic properties of large random matrices with independent entries. *J. Math. Phys.* **37** 5033–5060. MR1411619

[12] LINDEBERG, J. W. (1922). Eine neue herleitung des exponentialgesetzes in der wahrscheinlichkeitsrechnung. *Math. Z.* **15** 211–225. MR1544569

[13] MOSSEL, E., O'DONNELL, R. and OLESZKIEWICZ, K. (2005). Noise stability of functions with low influences: Invariance and optimality. Available at http://arxiv.org/math.PR/0503503.

[14] PASTUR, L. A. (1972). The spectrum of random matrices. *Teoret. Mat. Fiz.* **10** 102–112. MR0475502

[15] PAULAUSKAS, V. and RAČKAUSKAS, A. (1989). *Approximation Theory in the Central Limit Theorem. Exact Results in Banach Spaces.* Kluwer, Dordrecht. MR1015294

[16] ROTAR, V. I. (1979). Limit theorems for polylinear forms. *J. Multivariate Anal.* **9** 511–530. MR0556909

[17] SCHENKER, J. H. and SCHULZ-BALDES, H. (2005). Semicircle law and freeness for random matrices with symmetries or correlations. *Math. Res. Lett.* **12** 531–542. MR2155229

[18] SUIDAN, T. (2006). A remark on a theorem of Chatterjee and last passage percolation. *J. Phys. A* **39** 8977–8981. MR2240468

[19] TALAGRAND, M. (2003). *Spin Glasses: A Challenge for Mathematicians. Cavity and Mean Field Models.* Springer, Berlin. MR1993891

[20] TROTTER, H. F. (1959). Elementary proof of the central limit theorem. *Archiv der Mathem.* **10** 226–234. MR0108847

[21] WIGNER, E. P. (1958). On the distribution of the roots of certain symmetric matrices. *Ann. of Math. (2)* **67** 325–327. MR0095527

[22] WILKINSON, J. H. (1965). *The Algebraic Eigenvalue Problem.* Clarendon Press, Oxford. MR0184422

[23] ZOLOTAREV, V. M. (1977). Ideal metrics in the problem of approximating the distributions of sums of independent random variables. *Theory Probab. Appl.* **22** 449–465. MR0455066

DEPARTMENT OF STATISTICS
367 EVANS HALL #3860
UNIVERSITY OF CALIFORNIA, BERKELEY
BERKELEY, CALIFORNIA 94720
USA
E-MAIL: sourav@stat.berkeley.edu
URL: http://www.stat.berkeley.edu/~sourav